%% file: main.tex
\documentclass{aero}

\input{sections/preamble}

\ADsetup{
title    = {Chance constraints transcription and failure risk estimation for stochastic trajectory optimisation},
author   = {Thomas Caleb$^{1}$\cor{}, Roberto Armellin$^{2}$, and St\'ephanie Lizy-Destrez$^{1}$},
address  = {
	1. DCAS/SaCLaB, ISAE-SUPAERO, Toulouse 31055, France \sep
	2. Te P\=unaha \=Atea -- Space Institute, University of Auckland, Auckland 1010, New Zealand
},
email    = {Thomas Caleb, thomas.caleb@dycsyt.com},
abstract = {
Stochastic trajectory optimisation under uncertainty requires robust constraint satisfaction through chance constraints.
However, existing transcription methods remain limited to scalar constraints or highly specific structures while introducing substantial conservatism.
This work presents two general-purpose transcription methods for multi-dimensional Gaussian chance constraints for trajectory optimisation problems under uncertainty.
The spectral radius method extends existing methods to arbitrary multi-dimensional constraints with reduced conservatism. The refined first-order method achieves superior tightness with linear complexity.
In addition, a $d$-th order risk estimation methodology provides conservative failure probability estimates with limited conservatism in high dimensions in quadratic complexity.
Applied to an optimal control with uncertainties setting, the first-order transcription achieves near-optimal fuel consumption while maintaining the failure risk below the target.
The spectral radius method incurs approximately \SI{0.7}{\kilo\gram} additional fuel consumption due to excessive conservatism and a 51\% increase in computational time due to its cubic complexity.
High-dimensional tests show that the proposed risk estimation method provides accurate risk estimates, while previously developed methods exhibit exponential growth in conservatism with respect to constraint dimension.},
keywords = {Multi-dimensional Gaussian constraints \sep Stochastic trajectory optimisation \sep Transcription methods \sep Failure risk estimation}
}
\begin{document}

\maketitle

\input{sections/1_introduction}
\input{sections/2_background}
\input{sections/3_methodology.tex}
\input{sections/4_application}
\input{sections/5_conclusion.tex}
\input{sections/7_funding_acknoledgments}
\appendix
\input{sections/6_appendix}

\section*{References}
\bibliographystyle{abbrvnat}
\bibliography{references}

\end{document}

%% file: sections/preamble.tex
\usepackage[utf8]{inputenc}
\usepackage{textcomp}

\usepackage[version=4]{mhchem}
\usepackage{siunitx}
\usepackage{lipsum}

\usepackage[nonumberlist,acronym]{glossaries} 
\usepackage{graphicx,setspace,parskip,bm,amsmath,amssymb,mathrsfs} 

\usepackage{amsthm}
\usepackage{longtable,tabularx,multicol,multirow,threeparttablex} 
\usepackage{float} 
\usepackage{makecell} 
\usepackage{booktabs} 
\usepackage{upgreek}
\usepackage{footnpag,marginnote} 
\usepackage{booktabs} 
\usepackage[font=small]{caption} 
\usepackage{subcaption} 
\usepackage{comment}
\usepackage{microtype} 
\usepackage{makecell} 
\usepackage{setspace} 
\usepackage{pifont}
\usepackage{xcolor} 
\usepackage{colortbl} 
\usepackage{mathtools}
\usepackage{centernot}
\usepackage{pgfplots}
\pgfplotsset{compat=newest}
\pgfplotsset{plot coordinates/math parser=false}
\newlength\figureheight
\newlength\figurewidth
\usepackage{environ,tikz} 
\usetikzlibrary{mindmap, trees, calc}
\usetikzlibrary{positioning,fit}
\usetikzlibrary{through}
\usetikzlibrary{shapes.geometric, arrows}
\usetikzlibrary{automata,positioning}

\tikzset{block/.style={draw,thick,text width=2cm,minimum height=1cm,align=center},
         line/.style={-latex}
}
\usepackage{etoolbox}
\usepackage{enumitem}
\usepackage{algorithm}
\usepackage{algpseudocode}
\usepackage[numbers]{natbib}

\newtheorem{theorem}{Theorem}

\newtheorem{remark}{Remark}

\newtheorem{proposition}{Proposition}

\DeclareSIUnit{\AU}{AU}
\DeclareSIUnit{\day}{d}
\DeclareSIUnit{\deg}{deg}
\DeclareSIUnit{\year}{yr}

\newcommand{\Eq}[1]{Eq.~\eqref{#1}}

\newcommand{\Fig}[1]{Fig.~\ref{#1}}

\newcommand{\Tab}[1]{Table~\ref{#1}}

\newcommand{\Prop}[1]{Proposition~\ref{#1}}
\newcommand{\Rmk}[1]{Remark~\ref{#1}}
\newcommand{\Rem}[1]{Remark~\ref{#1}}

\newcommand{\Thm}[1]{Theorem~\ref{#1}}

\newcommand{\App}[1]{Appendix~\ref{#1}}
\newcommand{\Sec}[1]{Section~\ref{#1}}

\newcommand{\ie}{i.\,e.,~}

\newcommand{\eg}{e.\,g.,~}

\newacronym{ABM}{ABM}{Adams-Bashforth-Moulton}
\newacronym{ADS}{ADS}{automatic domain splitting}
\newacronym{AUL}{AUL}{augmented Lagrangian}
\newacronym{CDF}{CDF}{cumulative distribution function}
\newacronym{CR3BP}{CR3BP}{circular restricted three-body problem}
\newacronym{EMS}{EMS}{Earth—Moon system}
\newacronym{GMM}{GMM}{Gaussian Mixture Model}
\newacronym{DRO}{DRO}{distant retrograde orbit}
\newacronym{NRHO}{NRHO}{near-rectilinear halo orbit}
\newacronym{DROs}{DROs}{distant retrograde orbits}
\newacronym{NRHOs}{NRHOs}{near-rectilinear halo orbits}
\newacronym{DA}{DA}{differential algebra}
\newacronym{DACE}{DACE}{differential algebra core engine}
\newacronym{DADDy}{DADDy}{differential algebra-based differential dynamic programming}
\newacronym{DDP}{DDP}{differential dynamic programming}
\newacronym{DNC}{DNC}{did not converge}
\newacronym{DOPRI8}{DOPRI8}{Dormand-Prince 8th-order embedded Runge-Kutta method}
\newacronym{DOP853}{DOP853}{Dormand-Prince of order 8(5,3) embedded Runge-Kutta method}
\newacronym{EOM}{EoM}{equations of motion}
\newacronym{GEO}{GEO}{geostationary orbit}
\newacronym{GTO}{GTO}{geostationary transfer orbit}
\newacronym{IC}{IC}{initial conditions}
\newacronym{JPL}{JPL}{Jet Propulsion Laboratory}
\newacronym{MC}{MC}{Monte-Carlo}
\newacronym{MB}{MB}{megabytes}
\newacronym{LAM}{LAM}{likelihood agreement measure}
\newacronym{LEO}{LEO}{low Earth orbit}
\newacronym{MEO}{MEO}{medium Earth orbit}
\newacronym{NASA}{NASA}{National Aeronautics and Space Administration}
\newacronym{NEA}{NEA}{near-Earth asteroid}
\newacronym{NEO}{NEO}{near-Earth object}
\newacronym{NSG}{NSG}{non-spherical gravity}
\newacronym{ODE}{ODE}{ordinary differential equation}
\newacronym{PECE}{PECE}{predictor-corrector}
\newacronym{PCK}{PCK}{planetary constants kernel}
\newacronym{PDF}{PDF}{probability density function}
\newacronym{PDS}{PDS}{Planetary Data System}
\newacronym{PN}{PN}{projected Newton}
\newacronym{RAAN}{RAAN}{right ascension of the ascending node}
\newacronym{RHS}{RHS}{right-hand side}
\newacronym{RK}{RK}{Runge-Kutta}
\newacronym{RK78}{RK78}{order 7(8) Runge–Kutta–Fehlberg embedded method}
\newacronym{RSS}{RSS}{root sum square}
\newacronym{RTN}{RTN}{radial-tangential-normal}
\newacronym{RPF}{RPF}{roto-pulsating frame}
\newacronym{STT}{STT}{state transition tensor}
\newacronym{STTs}{STTs}{state transition tensors}
\newacronym{STM}{STM}{state-transition matrix}
\newacronym{TDA}{TDA}{Taylor Differential Algebra}
\newacronym{ToF}{ToF}{time-of-flight}
\newacronym{TPBVP}{TPBVP}{two-point boundary value problem}
\newacronym{UOM}{uom}{unit of measurement}
\newacronym{VSVO}{VSVO}{variable-step variable-order}
\newacronym{VV}{V\&V}{verification and validation}


%% file: sections/1_introduction.tex
\section{Introduction}
\label{sec:introduction}
Space exploration has significantly advanced over the decades, with missions increasingly using complex dynamical systems, such as those involving multiple planetary swing-bys \citep{BoutonnetEtAl_2024_DtJT} or the challenging dynamics of the \gls*{CR3BP} \citep{SmithEtAl_2020_TAPaOoNAtRHttM}. These missions require optimal trajectory design, which involves minimising an objective function (such as fuel consumption or \gls*{ToF}) while ensuring robustness against the sensitivity to measurement, control, and model errors.

In trajectory optimisation, robust planning has traditionally involved deterministic approaches that ignore stochastic effects, which can lead to suboptimal or infeasible solutions when adding uncertainties. However, certain stochastic solvers incorporate transcription methods that transform stochastic cost functions into deterministic ones. Some methods even added constraints directly into the cost function, making them soft constraints, \ie constraints that may be violated \citep{YiEtAl_2019_NCCVDDP, OzakiEtAl_2018_SDDPwUTfLTTD}. On the other hand, hard constraints, \ie constraints that cannot be violated, need to be implemented separately and cannot be enforced with a probability of \num{100}\% in most cases, \eg if uncertainties are Gaussian \citep{BenedikterEtAl_2022_CAtCCwAtSLTTO}.

Recent research has focused on stochastic solvers that directly address uncertainties through chance constraints.
These constraints are relaxed hard constraints that allow for a given failure risk. However, they are not directly tractable. Approaches such as reinforcement learning \citep{FedericiEtAl_2021_DLTfASGdPO,HoltEtAl_2021_OQLVRLwGS}, distribution approximation \citep{GharaviEtAl_2024_PECAfADiHS}, \gls*{MC} sampling \citep{BlackmoreEtAl_2010_APPCAoCCSPC,GeletuEtAl_2013_AaAoCCAtSOuU}, or transcription methods can be used to make these constraints manageable. These transcription methods reformulate optimisation problems to enforce chance constraints, ensuring failure probabilities remain within acceptable thresholds \citep{BlackmoreEtAl_2011_CCOPPwO,RidderhofEtAl_2020_CCCCfLTMFTO,OzakiEtAl_2020_TSOCfNCTOP}. However, existing transcription methods often make conservative assumptions, such as bounding covariances with their spectral radii, which can sometimes lead to overly cautious solutions and potentially reduce performance \citep{BenedikterEtAl_2022_CAtCCwAtSLTTO,MarmoZavoli_2024_CCMfCCoLTIM}.

Despite these advancements, significant gaps remain in the literature. Many existing transcription methods are limited by their reliance on specific assumptions about constraints \citep{RidderhofEtAl_2020_CCCCfLTMFTO} and their inability to generalise to multi-dimensional constraints \citep{BlackmoreEtAl_2011_CCOPPwO,NakkaChung_2023_TOoCCNSSfMPuU}, such as maintaining a bounded thrust magnitude and keeping a safe distance from a celestial body. Addressing these limitations could enable the design of scalable, efficient, and generalisable methods for stochastic trajectory optimisation.

In this work, we present two general multidimensional transcription methods for Gaussian chance constraints with low conservatism compared to existing approaches. 
Additionally, we introduce a $d$-th order risk estimation methodology and a conservatism metric for comparing novel and state-of-the-art methods.

\Sec{sec:background} reviews notations and existing transcription methods. \Sec{sec:first_order} introduces the spectral radius-based and first-order transcription methods. \Sec{sec:risk_estimation} presents the $d$-th order risk estimation approach and a conservatism metric for method comparison. \Sec{sec:application} demonstrates the application of these methods to a standard stochastic trajectory optimisation and a high-dimensional test case. Finally, \Sec{sec:conclusion} concludes this study.

%% file: sections/2_background.tex
\section{Background}
\label{sec:background}
This section presents the framework and the notations of chance constraints.

Uncertainties in the state vector $\bm{x}$ of size $N_x$ are often modelled as a multivariate Gaussian variable with mean $\Bar{\bm{x}}$ and covariance $\bm{\Sigma}_x$, denoted as $\bm{x}\sim\mathcal{N}\left(\Bar{\bm{x}},\bm{\Sigma}_x\right)$ \citep{Tong_1990_TMND}. If the control $\bm{u}$ of size $N_u$ has a linear feedback gain $\bm{K}$, such that $\bm{u}=\Bar{\bm{u}} + \bm{K}(\bm{x}-\Bar{\bm{x}})$, then $\bm{u}\sim\mathcal{N}\left(\Bar{\bm{u}},\bm{\Sigma}_u\right)$, where $\Bar{\bm{u}}$ is the nominal control and $\bm{\Sigma}_u=\bm{K}\bm{\Sigma}_x\bm{K}^{\textrm{T}}$.
As a result, enforcing hard constraints may become infeasible in all cases because $\bm{x}$ and $\bm{u}$ are unbounded, as is:
\begin{equation}
    \label{eq:hard_constraint}
    \|\bm{u}\| \leq u_{\max}.
\end{equation}
The chance constraints paradigm allows us to enforce this type of constraint.
Given an acceptable failure risk $\beta\in]0, 1[$, the constraint in \Eq{eq:hard_constraint} becomes:
\begin{equation}
    \label{eq:chance_constraint}
    \mathbb{P}\left(\|\bm{u}\| \leq u_{\max}\right) \geq 1-\beta.
\end{equation}
The probability that the control satisfies \Eq{eq:hard_constraint}, \ie the probability that the control is feasible, must be higher than $1-\beta$.

However, enforcing \Eq{eq:hard_constraint} while solving an optimisation problem is intractable in practice and computing the derivatives of $\mathbb{P}\left(\|\bm{u}\| \leq u_{\max}\right)$ with respect to the states and control variables is cumbersome. 
To make the problem tractable, transcription methods are introduced to transform chance constraints into tractable deterministic ones.  
Many stochastic optimisation techniques rely on such transcriptions to approximate chance constraints of the form $\mathbb{P}\left(\bm{y} \preceq \bm{0}\right) \geq 1 - \beta$ with tractable deterministic sufficient conditions $\mathcal{C}\left(\Bar{\bm{y}}, \bm{\Sigma}_{\bm{y}}, \beta\right)$. We then have: $\mathcal{C}\left(\Bar{\bm{y}}, \bm{\Sigma}_{\bm{y}}, \beta\right) \implies \mathbb{P}\left(\bm{y} \preceq \bm{0}\right) \geq 1 - \beta$.
Since $\mathcal{C}$ is only a sufficient condition, it introduces conservatism, i.e., additional safety margins to ensure constraint satisfaction. Several transcription methods already exist in the literature. Nevertheless, none of them directly address multidimensional chance constraints in their full generality.

\subsection{Tailored norm-constraint transcriptions}
\citet{RidderhofEtAl_2020_CCCCfLTMFTO} introduce sufficient conditions specifically designed for chance constraints of the form $\mathbb{P}\left(\|\bm{u}\| \leq u_{\max}\right) \geq 1 - \beta$, where $\bm{u}$ is a Gaussian random variable, $\|\cdot\|$ is the Euclidean norm, and $u_{\max}$ denotes a hard bound on control magnitude. Such constraints are common in space trajectory optimisation \citep{OzakiEtAl_2020_TSOCfNCTOP,BenedikterEtAl_2022_CAtCCwAtSLTTO,BooneMcMahon_2022_NGCCTCUGMaRA}.
\citet{OguriLantoine_2022_SSCPfRLTTDuU} propose a refined version, replacing the conservative margin in \citet{RidderhofEtAl_2020_CCCCfLTMFTO} with a tighter one based on the quantile of the chi-squared distribution. They also suggest estimating the spectral radius using the Frobenius norm of the Cholesky decomposition \citep{Cholesky_1910_SLRNDSDL}, though this does not reduce the overall complexity.
    
These methods require no structural assumptions on the distribution and rely solely on the mean and covariance of $\bm{u}$. They involve computing the spectral radius of the covariance matrix, with computational complexity $\mathcal{O}(N_u^3)$ when using Jacobi eigenvalue algorithm.
Despite their conservatism, these tailored norm-constraint transcriptions are often favoured for their simplicity and minimal assumptions, and have been adopted in several studies on stochastic optimisation in astrodynamics \citep{BenedikterEtAl_2022_CAtCCwAtSLTTO,MarmoZavoli_2024_CCMfCCoLTIM}.

\subsection{Linear one-dimensional transcriptions}
A general approach for one-dimensional linear chance constraints involving multivariate Gaussian variables was introduced by \citet{BlackmoreEtAl_2011_CCOPPwO}. For constraints of the form $\mathbb{P}\left(\bm{h}^{\textrm{T}} \bm{z} \leq a \right) \geq 1 - \beta$, with $\bm{z} \sim \mathcal{N}(\Bar{\bm{z}}, \bm{\Sigma}_{\bm{z}})$ and $\bm{h} \in \mathbb{R}^{N_z}$, $a \in \mathbb{R}$, the method provides an exact deterministic transcription using the Gaussian quantile function. It has been widely used in stochastic control and trajectory optimisation \citep{FarinaEtAl_2016_SLMPCwCC–aR,OkamotoEtAl_2018_OCCfSSuCC,RidderhofEtAl_2019_NUCwICS,OzakiEtAl_2020_TSOCfNCTOP,LewEtAl_2020_CCSCPfRTO,OguriMcMahon_2021_RSGaSBuUSOCA,BooneMcMahon_2022_NGCCTCUGMaRA}.
\citet{NakkaChung_2023_TOoCCNSSfMPuU} propose a similar, but more conservative condition, replacing the Gaussian quantile with a bound that avoids the numerical evaluation of the inverse Gaussian \gls*{CDF}.
    
These methods require either linear constraints or local linearisability, which assumes the covariance is sufficiently small. Under this assumption, norm-type constraints can also be transcribed using this approach, by setting $\bm{h} = -\frac{\Bar{\bm{u}}}{\|\Bar{\bm{u}}\|}$, $a = u_{\max} - \|\Bar{\bm{u}}\|$, and $\bm{z} = \delta\bm{u} \sim \mathcal{N}(\bm{0}, \bm{\Sigma}_{\bm{u}})$. The computational complexity is dominated by matrix-vector multiplications, scaling as $\mathcal{O}(N_u^2)$.

%% file: sections/3_methodology.tex
\section{{Generic $d$\textendash{dimensional} transcription methods}}
\label{sec:first_order}
This section presents the development of two $d$\textendash{dimensional} general transcription methods. In contrast, the state-of-the-art is only able to tackle either specific constraint types on normal multidimensional distributions or general linear single-dimensional constraints.
Let us define $\mathbb{N}_d=\{1,2,\dotsc, d\}$ and for $\left(\bm{y},\bm{z}\right)\in\left(\mathbb{R}^d\right)^2$, we define $\bm{y} \preceq \bm{z} \Leftrightarrow \forall i\in\mathbb{N}_d, \ y_i \leq z_i$.

The goal is to transcribe a general multi-dimensional inequality constraint such as: $\mathbb{P}\left(\bm{g}\left(\bm{x}, \bm{u}\right) \preceq \bm{0} \right) \geq 1-\beta$, where $\bm{g}:\mathbb{R}^{N_x}\times\mathbb{R}^{N_u}\longrightarrow\mathbb{R}^{d}$ is a differentiable function. We assume that $\bm{g}$ can be linearised with respect to the state $\bm{x}$ and control $\bm{u}$. This is exact if $\bm{g}$ is linear, and otherwise justified when the uncertainty is sufficiently small so that a local linear approximation is accurate. This assumption, referred to as the weakly non-linear assumption, is adopted throughout the remainder of this work.
We define $\bm{y}=\bm{g}(\bm{x},\bm{u})$. 
If $\delta\bm{x}=\bm{x}-\Bar{\bm{x}}$, then:
\begin{equation}
    \bm{y} = \bm{g}\left(\Bar{\bm{x}}, \Bar{\bm{u}}\right) + \left[\nabla _{x}\bm{g} + \nabla _{u}\bm{g}\bm{K}\right]\delta\bm{x},
\end{equation}
where $\nabla \bm{g}$ is the gradient of $\bm{g}$, $\nabla_x\bm{g}$ is the matrix containing the $N_x$ first columns of $\nabla \bm{g}$, and $\nabla_u \bm{g}$ is the matrix containing its $N_u$ last columns. Thus, the random variable $\bm{y}$ is such that $\bm{y}\sim\mathcal{N}(\Bar{\bm{y}}, \bm{\Sigma}_{y})$ with:
\begin{equation}
    \begin{aligned}
        \Bar{\bm{y}} &= \bm{g}(\Bar{\bm{x}}, \Bar{\bm{u}}), \\
        \bm{\Sigma}_{y} &= \left(\nabla_x \bm{g}+\nabla_u \bm{g}\bm{K}\right)
        \bm{\Sigma}_x \left(\nabla_x \bm{g}+\nabla_u \bm{g}\bm{K}\right)^{\textrm{T}}.
    \end{aligned}
\end{equation}
Indeed, note that $\delta\bm{x}\sim\mathcal{N}(\bm{0}, \bm{\Sigma}_{x})$, therefore, since $\bm{y}$ is a linear combination of a Gaussian variable, it is Gaussian with the given mean and covariance. 

Then, the chance constraint $\mathbb{P}\left(\bm{y} \preceq \bm{0} \right) \geq 1-\beta $ is covered by \citet{BlackmoreEtAl_2011_CCOPPwO} if $d=1$. The idea of the following results is to generalise this result to any dimension $d$.
Let $\rho(\bm{M})$ be the positive square root of the spectral radius of the matrix $\bm{M}$ and $\Phi^{-1}_d$ be the inverse \gls*{CDF} of chi-squared distribution with $d$ degrees of freedom \citep{AbramowitzStegun_1964_HoMFwFGaMT}. We also define: $\Psi^{-1}_d\left(\beta\right) = \sqrt{\Phi^{-1}_d\left(1-\beta\right)}$, and for $r\geq0$: $\Psi_d\left(R\right) = 1 - \Phi_d\left(R^2\right)$.
\begin{theorem}[\textsc{Spectral radius transcription.}]
\label{thm:spectral_radius_transcription}
\begin{equation}
    \label{eq:thm:spectral_radius_transcription}
        \forall i \in\mathbb{N}_d, \ \Bar{y}_i  +  \Psi_d^{-1}\left(\beta\right)\rho\left(\bm{\Sigma}_{y}\right)\leq 0  \implies \mathbb{P}\left(\bm{y}\preceq \bm{0}\right) \geq 1-\beta
\end{equation}
\end{theorem}
\begin{proof}
    All the proofs of the results presented in this work are reported in \App{sec:appendix}
\end{proof}
This method has cubic computational complexity in $d$, due to the eigenvalue decomposition.
\begin{remark}
\label{rmk:oguri}
Under the weakly non-linear assumption, the spectral radius transcription of \Thm{thm:spectral_radius_transcription} is less conservative than that of \citet{OguriLantoine_2022_SSCPfRLTTDuU} when applied to Gaussian norm constraints, and therefore also less conservative than those of \citet{RidderhofEtAl_2020_CCCCfLTMFTO}.
\end{remark}

In the remainder of this work, it will be shown that this sufficient condition is highly conservative. This motivates the development of less conservative transcription methods capable of handling multi-dimensional chance constraints. Let $\bm{\sigma}_{y}$ be the vector of the positive square roots of the diagonal coefficients of $\bm{\Sigma}_{y}$.
The authors present so-called first-order sufficient conditions to satisfy the chance constraints:
\begin{theorem}[\textsc{First-order transcription.}]
    \label{thm:first_order_transcription}
    \begin{equation}
        \label{eq:first_order_theorem}
        \Bar{\bm{y}} + \Psi_d^{-1}\left(\beta\right)\bm{\sigma}_{y} \leq \bm{0} 
        \implies \mathbb{P}\left(\bm{y}\preceq \bm{0}\right) \geq 1-\beta
    \end{equation}
\end{theorem}
Compared to the spectral radius method, this method only requires computing the square roots of $d$ terms of a known matrix, therefore, the complexity is $\mathcal{O}(d)$.
\begin{remark}
    \label{rmk:spectral_vs_first_order}
    The transcription in \Thm{thm:first_order_transcription} is less conservative than that in \Thm{thm:spectral_radius_transcription}. For $d = 1$, these transcriptions are equivalent.
\end{remark}

\section{{Risk estimation methods and conservatism metric}}
\label{sec:risk_estimation}
\subsection{Risk estimation}
While deterministic transcriptions provide conservative guarantees, it is often desirable to estimate the actual failure risk $\beta_\mathrm{R} = 1 - \mathbb{P}(\bm{y} \preceq \bm{0})$ more accurately. This section introduces risk estimation techniques that approximate the true probability of Gaussian chance constraint violation.
Although the risk can, in principle, be estimated using numerical integration or Monte Carlo methods \citep{RobertCasella_2004_MCSM}, these approaches are generally intractable during the optimisation process.

To address this, we propose a methodology to conservatively estimate the failure risk, \ie to provide an upper bound of $\beta_{\textrm{R}}$.
This allows us to associate a risk estimate with each transcription method by identifying (when possible) the smallest $\beta_{\textrm{T}}$ such that the sufficient condition $\mathcal{C}\left(\Bar{\bm{y}}, \bm{\Sigma}_{\bm{y}}, \beta_{\textrm{T}}\right)$ holds.
Since transcriptions correspond to sufficient conditions, we always have $\beta_{\textrm{T}} \geq \beta_{\textrm{R}}$, ensuring that the true risk is not underestimated. Let $\Phi_G$ be the \gls*{CDF} of the standard Gaussian distribution and $\bm{r} = \left[-\frac{\Bar{y}_1}{\sigma_{\bm{y},1}}, \dotsc, -\frac{\Bar{y}_d}{\sigma_{\bm{y},d}}\right]^{\mathrm{T}}$.
\begin{proposition}
\label{prop:risk_estimates_all}
Let us denote $\delta\|\bm{u}\|=\dfrac{\|\Bar{\bm{u}}\| - u_{\max}}{\rho(\bm{\Sigma}_{\bm{u}})}$ for conciseness. The closed-form expressions in \Tab{tab:risk_estimations} yield analytical estimates $\beta_{\textrm{T}}$ that are upper bounds on the true risk $\beta_{\textrm{R}} = 1 - \mathbb{P}(\bm{y} \preceq \bm{0})$, under the assumptions specific to each method and using the notations given in the corresponding references.
\begin{table}[htbp]
    \centering
    \caption{Failure risk estimates for various constraint transcriptions.}
    \label{tab:risk_estimates}
    
    \begin{tabular}{ l c c c p{5cm} p{5cm} }
    \toprule
    \textbf{Transcription method} &  \textbf{Symbol} & \textbf{Condition} & \textbf{Risk Estimate Expression} \\
    \midrule
    \multirow{2}{*}{\citet{RidderhofEtAl_2020_CCCCfLTMFTO}} & \multirow{2}{*}{$\beta_{\textrm{T},\textrm{R}}$} &
    $ \|\Bar{\bm{u}}\| \leq u_{\max}$ and
    $N_u \in \{1,2\}$ &
    $\exp\left[-\dfrac{1}{2}\left(\delta\|\bm{u}\|\right)^2\right]$ \\
    &  &
    $\|\Bar{\bm{u}}\| \leq u_{\max}$ and
    $N_u > 2$ and $\left(\delta\|\bm{u}\| + \sqrt{N_u}\right) \leq 0$ &
    $\exp\left[-\dfrac{1}{2}\left(\delta\|\bm{u}\| + \sqrt{N_u}\right)^2\right]$ \\
    \midrule
    \citet{OguriLantoine_2022_SSCPfRLTTDuU} & $\beta_{\textrm{T},\textrm{OL}}$ &
    $\|\Bar{\bm{u}}\| \leq u_{\max}$ &
    $\Psi_{N_u}\left(\delta\|\bm{u}\|\right)$ \\
    \midrule
    \citet{NakkaChung_2023_TOoCCNSSfMPuU} & $\beta_{\textrm{T},\textrm{NC}}$ &
    Linearity and $\bm{h}^{\textrm{T}}\Bar{\bm{z}} \leq a$ &
    $\dfrac{\bm{h}^{\textrm{T}}\bm{\Sigma}_{\bm{z}}\bm{h}}{\bm{h}^{\textrm{T}}\bm{\Sigma}_{\bm{z}}\bm{h} + (a - \bm{h}^{\textrm{T}}\Bar{\bm{z}})^2}$ \\
    \midrule
    \citet{BlackmoreEtAl_2011_CCOPPwO} & $\beta_{\textrm{T},\textrm{B}}$ &
    Linearity &
    $1 - \Phi_G\left(\dfrac{a - \bm{h}^{\textrm{T}}\Bar{\bm{z}}}{\sqrt{\bm{h}^{\textrm{T}}\bm{\Sigma}_{\bm{z}}\bm{h}}}\right)$ \\
    \midrule
    Spectral Radius (\Thm{thm:spectral_radius_transcription}) & $\beta_{\textrm{T},\rho}$ &
    Linearity and $\Bar{\bm{y}} \preceq \bm{0}$ &
    $\Psi_d\left[\dfrac{\min(-\Bar{\bm{y}})}{\rho(\bm{\Sigma}_{\bm{y}})}\right]$ \\
    \midrule
    First-Order (\Thm{thm:first_order_transcription}) & $\beta_{\textrm{T},1}$ &
    Linearity and $\Bar{\bm{y}} \preceq \bm{0}$ &
    $\Psi_d\left(\min\bm{r}\right)$ \\
    \bottomrule
    \end{tabular}
    \label{tab:risk_estimations}
\end{table}
\end{proposition}

We now introduce the so-called $d$\textendash{th}-order low-conservatism risk estimate, which refines the first-order approximation by incorporating higher-order effects.
Let us define $\Tilde{r}_1 \leq \cdots \leq \Tilde{r}_d$, the ordered components of the vector $\bm{r}$.
\begin{theorem}[\textsc{$d$\textendash{th}-order failure risk estimation}]
\label{thm:dth_order_risk}
If $\Bar{\bm{y}} \preceq \bm{0}$, then:
\begin{equation}
    \label{eq:def_beta_d}
    \beta_\mathrm{T} = 1 - \sum_{i=1}^d \left[\Psi_d(\Tilde{r}_i) - \Psi_d(\Tilde{r}_{i-1})\right] \max\left[0, 1 - \frac{1}{2} \sum_{j=1}^{i-1} I\left(\frac{d-1}{2}, \frac{1}{2}; \frac{\Tilde{r}_j}{\Tilde{r}_i}\right)\right] \geq \beta_\mathrm{R},
\end{equation}
where $I(a,b;x)$ denotes the regularised incomplete beta function: $I\left(a, b;x\right) = \frac{1}{B\left(a, b\right)} \int_{0}^{x} {t^{a-1}(1-t)^{b-1}dt}$, with $B(a,b)$ the Euler beta function \citep{AbramowitzStegun_1964_HoMFwFGaMT}.
\end{theorem}
Computing $\beta_{\textrm{T},d}$ requires at most $\mathcal{O}(d^2)$ operations, therefore, this transcription method has a complexity $\mathcal{O}(d^2)$.

\subsection{Conservatism metric}
To quantify the conservatism introduced by a transcription method, we define the metric $\gamma(\beta_\mathrm{T})$ as follows, comparing the estimated and the true failure risks: $\gamma(\beta_\mathrm{T}) = \frac{\beta_\mathrm{T}}{\beta_\mathrm{R}}  \sqrt{\frac{1 - \beta_\mathrm{R}^2}{1 - \beta_\mathrm{T}^2}}$.
The function $\gamma$ is monotonically increasing with respect to $\beta_\mathrm{T}$. Since by design any transcription method satisfies $\beta_\mathrm{T} \geq \beta_\mathrm{R}$, it follows that $\gamma(\beta_\mathrm{T}) \geq 1$.
Transcription methods can be compared using this metric: a method is said to be less conservative if its associated $\gamma$ value is closer to 1. When $1 - \beta_\mathrm{T} \ll 1$ while $\beta_\mathrm{R} \ll 1$, the transcription introduces excessive margins despite the actual failure risk being low, and is therefore penalised by a high conservatism value.
Moreover, when $\beta_\mathrm{T} \ll 1$ and $\beta_\mathrm{R} \ll 1$, we have the asymptotic approximation: $\gamma\left(\beta_{\textrm{T}}\right)\sim\frac{{\beta_{\textrm{T}}}}{\beta_{\textrm{R}}}$. 
Thus, $\gamma$ accurately captures the overestimation of the risk and serves as a relevant metric to compare different transcription methods.
We denote by $\gamma_{\rho}$, $\gamma_1$, and $\gamma_d$ the conservatism values of the spectral radius, first-order, and $d$th-order transcription methods, respectively.
\begin{theorem}{\textsc{Transcription methods hierarchy.}}
    \label{thm:conservatism_hierarchy}
    If $\Bar{\bm{y}}\preceq\bm{0}$, then: $1\leq\gamma_d\leq\gamma_1\leq\gamma_{\rho}$. Moreover, if $d = 1$, then: $\gamma_d=\gamma_1=\gamma_{\rho}$.
\end{theorem}
\Thm{thm:conservatism_hierarchy} formally establishes the superiority of the $d$th-order transcription method over both the first-order and the spectral radius approaches in terms of conservatism.

%% file: sections/4_application.tex
\section{Numerical application}
\label{sec:application}

This section presents numerical applications of the proposed methods. First, the spectral radius and first-order transcriptions are compared with existing one-dimensional methods from the literature in a chance-constrained setting. 
Then, the two transcriptions are compared within a full stochastic optimal control problem involving two-dimensional constraints, which are not addressed by existing methods in the literature.
Finally, the fidelity of the risk estimation method is evaluated on a high-dimensional test case.
\subsection{Control-norm chance constraints}
All the transcription methods are now compared on a control-norm chance constraint defined as follows:   $u_{\max} = \text{\SI{0.5}{\newton}}$, $\Bar{\bm{u}} = [\text{\num{0.3}}, \text{\num{0.37}}, \text{\num{-0.15}}]^{\textrm{T}} \text{ \SI{}{\newton}}$, and $\bm{\Sigma}_u = 0.1\mathbb{I}_3 + 10^{-3}(\bm{1}_{3}-\mathbb{I}_3) \text{\SI{}{\milli\newton\squared}}$, where $\mathbb{I}_3$ is the identity of size 3 and $\bm{1}_{3}$ is the all-ones matrix.
The nominal value of the control norm was set to $\|\Bar{\bm{u}}\| = u_{\max} -$ \SI{0.60}{\milli\newton}, so that the constraint is nearly saturated for illustrative purposes. Note that the magnitude of the deviations is on the order of $\sqrt{\|\bm{\Sigma}_u \|_2}$, \ie less than \SI{1}{\milli\newton}, thus satisfying the weakly non-linear hypothesis.

\begin{figure}[h]
    \centering
    \includegraphics[width=0.5\linewidth]{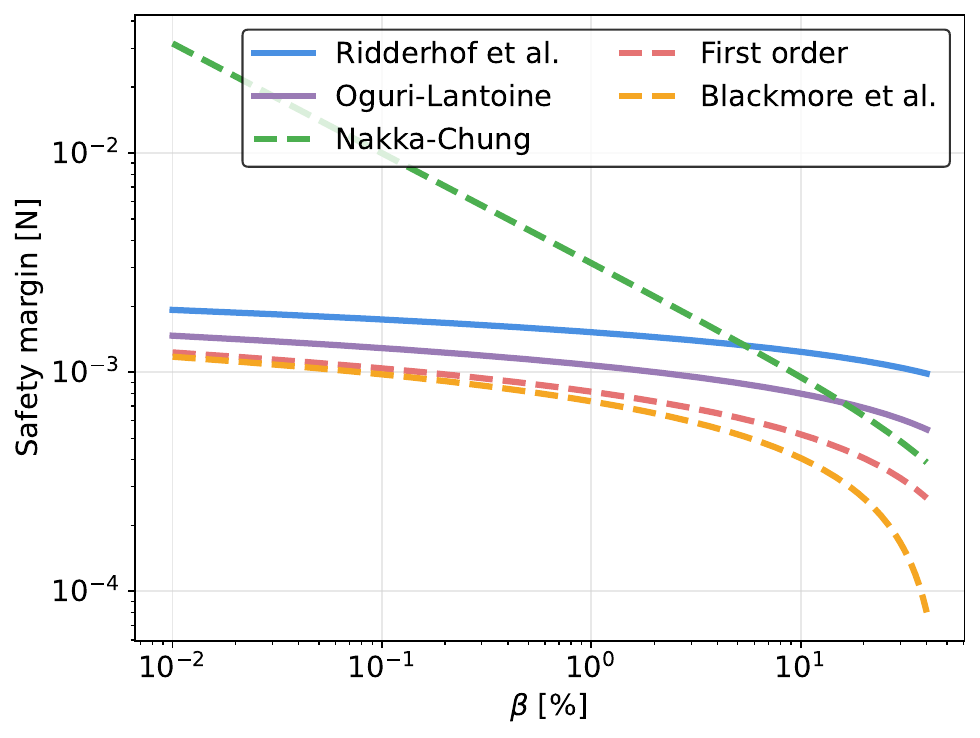}
    \caption{Values of the safety margins of the presented transcriptions.}
    \label{fig:safety_margins}
\end{figure}

We are interested in comparing the safety margins $\Delta \|\bm{u}\|(\Bar{\bm{u}},\bm{\Sigma}_u ,\beta)$ required by each transcription method: $\|\Bar{\bm{u}}\| + \Delta \|\bm{u}\| \left(\Bar{\bm{u}},\bm{\Sigma}_u ,\beta\right) \geq u_{\max} \implies \mathbb{P}\left( \|\bm{u}\|\leq u_{\max}\right) \geq 1-\beta$.
\Fig{fig:safety_margins} shows the safety margins obtained using the following methods: the tailored norm-constraint from \citet{RidderhofEtAl_2020_CCCCfLTMFTO} in blue, the tailored norm-constraint from \citet{OguriLantoine_2022_SSCPfRLTTDuU} in purple, the result of \citet{NakkaChung_2023_TOoCCNSSfMPuU} in green, the proposed \Thm{thm:spectral_radius_transcription} and \Thm{thm:first_order_transcription} in red, and the necessary and sufficient margin of \citet{BlackmoreEtAl_2011_CCOPPwO} in gold. The dashed lines highlight transcription methods that require the weakly non-linear.
This figure confirms that the transcription of \citet{OguriLantoine_2022_SSCPfRLTTDuU} is sharper than that of \citet{RidderhofEtAl_2019_NUCwICS}, as expected. Similarly, as shown in \Rmk{rmk:oguri}, the first-order and spectral radius transcriptions (which coincide for $d = 1$) are sharper than that of \citet{OguriLantoine_2022_SSCPfRLTTDuU}. The largest margin, for small values of $\beta$, is given by \citet{NakkaChung_2023_TOoCCNSSfMPuU}, which also requires linearisation.
The transcription of \citet{BlackmoreEtAl_2011_CCOPPwO}, under the small covariance assumption, yields the smallest safety margin that guarantees $\mathbb{P}( \|\bm{u}\| \leq u_{\max} ) \geq 1 - \beta$. However, note that the first-order transcription is the closest to this theoretical optimum.

We now estimate the failure risk, \ie $\beta_{\textrm{R}} = \mathbb{P}( \|\bm{u}\| > u_{\max} )$, using each method. \Tab{tab:transcription_test_case_control} compares the methods of \Prop{prop:risk_estimates_all}.
The value of $\beta_{\textrm{R}}$ is estimated using a \gls*{MC} sample of size \num{e7}.
\begin{table}[h]
    \centering
    \caption{Failure risk estimation comparison for the control-norm chance constraint.}
    \small
    \begin{tabular}{l c c c }
        \toprule
        \textbf{Type} & \textbf{Estimation method} & \textbf{Failure risk [\%]} & \textbf{$\gamma$} \\
        \midrule
        Numerical & $\beta_{\textrm{R}}$ (\gls*{MC}) & \num{2.89} & \textendash \\ 
        \midrule
        \multirow{2}{*}{\parbox{1cm}{Tailored}} & $\beta_{\textrm{T},\textrm{R}}$ & \num{98.9} & \num{232.6} \\
        & $\beta_{\textrm{T},\textrm{OL}}$ & \num{31.6} & \num{11.53} \\
        \midrule
        \multirow{3}{*}{\parbox{1cm}{Linear}} & $\beta_{\textrm{T},\textrm{NC}}$ & \num{21.7} & \num{7.695} \\
        &$\beta_{\textrm{T},1}$ & \num{5.77} & \num{1.999} \\ 
        &$\beta_{\textrm{T},\textrm{B}}$ & \num{2.89} & \num{0.9981} \\
        \bottomrule
    \end{tabular}
    \label{tab:transcription_test_case_control}
\end{table}
The results exhibit the same hierarchy as in \Fig{fig:safety_margins}. Specifically, the tailored norm-constraint methods, \ie from \citet{RidderhofEtAl_2020_CCCCfLTMFTO} and \citet{OguriLantoine_2022_SSCPfRLTTDuU}, yield highly conservative estimates due to the weak assumptions required for their applicability.
The method of \citet{NakkaChung_2023_TOoCCNSSfMPuU} provides a more accurate estimate than the tailored norm-constraint methods, consistent with the real safety margin being \SI{0.60}{\milli\newton}. As seen in \Fig{fig:safety_margins}, the associated failure risk is lower than those of \citet{RidderhofEtAl_2020_CCCCfLTMFTO} and \citet{OguriLantoine_2022_SSCPfRLTTDuU}.
The first-order risk estimate is less conservative than all others, except for the theoretical lower bound given thanks to \citet{BlackmoreEtAl_2011_CCOPPwO}. The latter slightly underestimates $\beta_{\textrm{R}}$ compared to the large-sample \gls*{MC} estimate, yielding a conservatism slightly under \num{1}, highlighting that linearised methods do not yield exact results.
Overall, this table shows that all transcription methods in the literature, except the estimate from \citet{BlackmoreEtAl_2011_CCOPPwO}, are more conservative than the first-order method.
However, even though the estimate from \citet{BlackmoreEtAl_2011_CCOPPwO} is superior, it cannot be applied to general multi-dimensional Gaussian chance constraints.

\subsection{Weakly non-linear stochastic optimal control problem}
This numerical application consists of computing a fuel-optimal Earth--mars transfer using the \gls*{DADDy} solver, a \gls*{DDP} solver presented in \citet{CalebEtAl_2025_TPBCSfFOLTTO}. 
The Earth--Mars problem is similar to the one solved in \citet{LantoineRussell_2012_AHDDPAfCOCPP2A}. Yet, the initial state is considered as a Gaussian variable $\bm{x}_{0}$ of covariance $\bm{\Sigma}_{\bm{x},0}=\textrm{diag}\left(10^{-6},10^{-6},10^{-6},5\times10^{-7},5\times10^{-7},5\times10^{-7}\right)$. In addition an unbiased state measurement noise is added at each time step with a covariance $\bm{\mathcal{Q}}_{\bm{x},k} = \bm{\Sigma}_{\bm{x},0}/10^4$. These values for the covariance matrices are small to ensure the weakly non-linear hypothesis is satisfied. The chance-constraints satisfied at each time step are:
\begin{equation}
    \mathbb{P}\left[\bm{g}(\bm{x},\bm{u}) = \left[\bm{u}^{\textrm{T}}\bm{u} - u_{\textrm{max}}^2, m_{\textrm{dry}} - m \right]^{\textrm{T}} \preceq\bm{0}\right]\geq 1-\beta,
\end{equation}
\ie the control magnitude should not exceed $u_{\max}$ and the total spacecraft mass $m$ should not be inferior to the spacecraft's dry mass. Note that no existing chance-constraints transcription can tackle this two-dimensional constraint.

\begin{figure}
    \centering
    \begin{subfigure}{.5\linewidth}
        \centering
        \includegraphics[width=\linewidth]{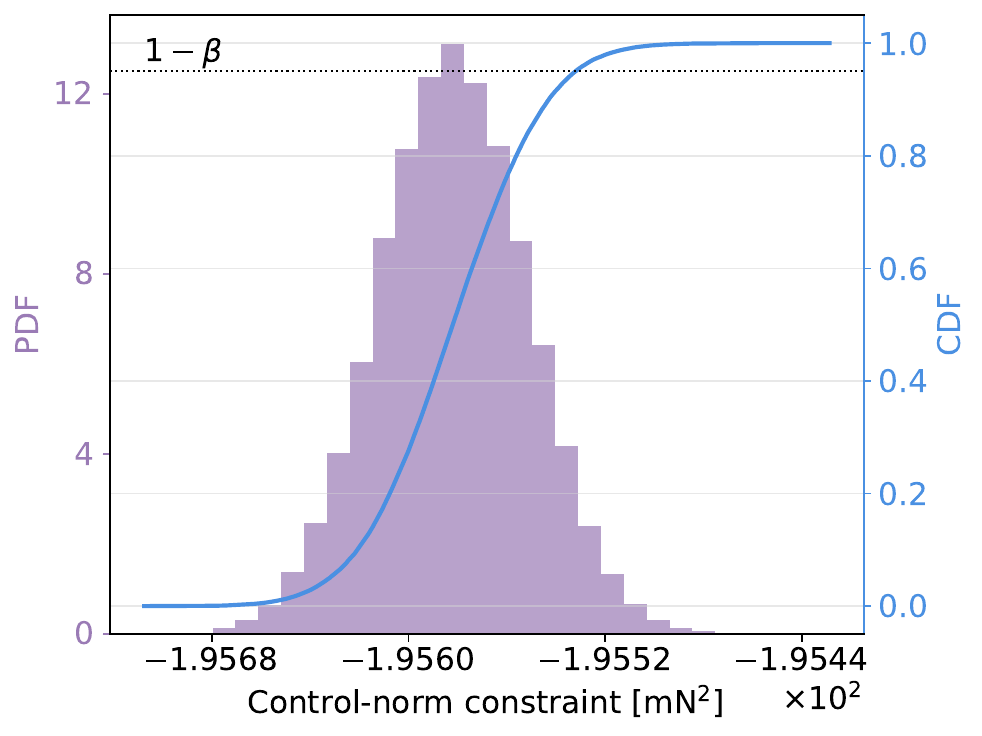}
        \caption{Spectral radius transcription.}
            \label{fig:transcription_earth_mars_distribution_sr}
    \end{subfigure}%
    \begin{subfigure}{.5\linewidth}
        \centering
        \includegraphics[width=\linewidth]{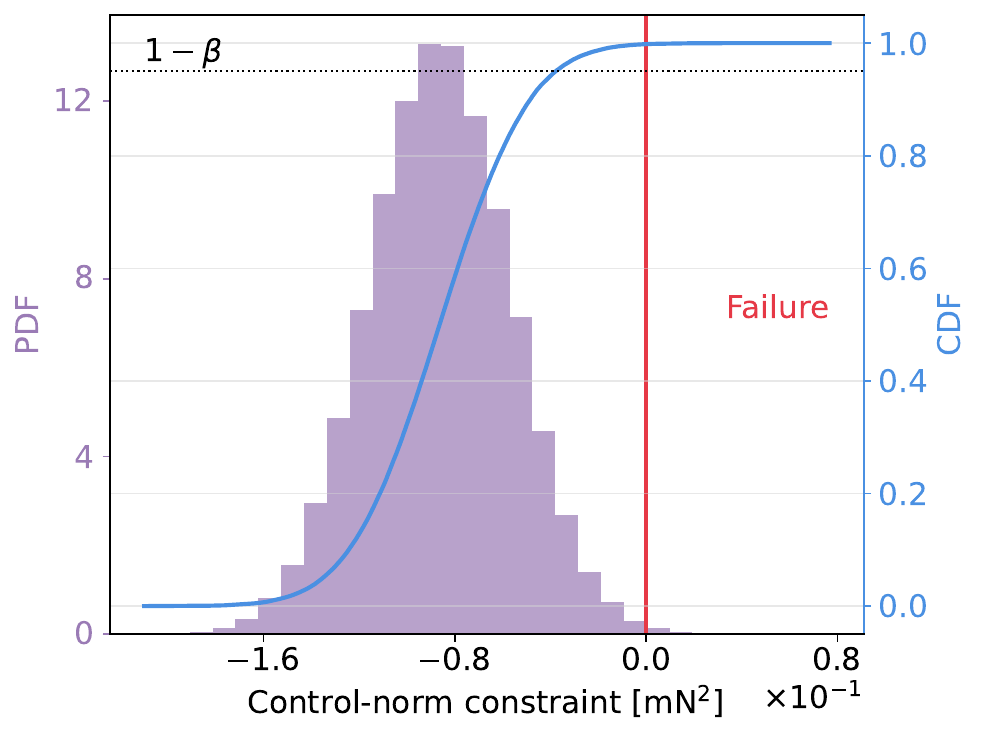}
        \caption{First-order transcription.}
            \label{fig:transcription_earth_mars_distribution_fo}
    \end{subfigure}
    \caption{Empirical \gls*{PDF} and \gls*{CDF} of the control-norm constraint violation for the two different transcriptions.}
    \label{fig:transcription_earth_mars_distribution}
\end{figure}

This optimisation problem is solved using the \gls*{DADDy} solver, for a target failure risk $\beta = 5\%$, and transcribed using both the spectral radius and the first-order methods.
The empirical \gls*{PDF} and the \gls*{CDF} of the control-norm constraint violations obtained after solving the stochastic optimal control problem are presented in \Fig{fig:transcription_earth_mars_distribution} and were obtained using Monte-Carlo sampling with \num{e5} realisations \citep{RobertCasella_2004_MCSM}.
\Fig{fig:transcription_earth_mars_distribution_sr} shows the results obtained using the spectral radius transcription, whereas \Fig{fig:transcription_earth_mars_distribution_fo} corresponds to the first-order transcription.
These qualitative results show that both distributions have similar standard deviations, namely approximately \SI{e-1}{\milli\newton\squared}. 
However, a non-negligible portion of the first-order transcription distribution lies in the failure region ($\bm{y} \npreceq \bm{0}$), while still satisfying the failure probability constraint ($\beta_{\mathrm{R}}$ remains below the target $\beta$). This behaviour is illustrated by the \gls*{CDF} intersecting the horizontal line $y = 1 - \beta$ before entering the failure region.
Conversely, the spectral radius transcription yields a distribution that is several orders of magnitude farther from the failure region, and is therefore entirely feasible. This indicates a high level of conservatism, as the available margin in the constraints is not exploited by the optimisation solver. This suggests that the spectral radius transcription may not accurately propagate these margins within the optimisation process.

Quantitative results comparing the two approaches are shown in \Tab{tab:earth_mars_L-SODA_transcription}. 
\begin{table}[h]
    \centering
    \caption{Convergence data for the stochastic Earth--Mars transfer ($\beta = 5\%$).}
    \begin{tabular}{l c c c c c c c} 
        \toprule
        \textbf{Transcription} & \textbf{RT [\SI{}{\second}]} & \textbf{\gls*{DDP} $n$}  & \textbf{Nominal error} & \textbf{Mean cost [\SI{}{\kilo\gram}]} & \textbf{$\beta_{\textrm{R}}$ [\%]} & \textbf{$\gamma$} \\
        \midrule
        Spectral radius & \num{69.2} & \num{371} & \num{-6e-6} & \num{397.6} & $<$\num{3e-3} & $>$\num{1700} \\
        First order & \num{45.9} & \num{370} & \num{-1e-10} & \num{396.9} & \num{2.93} & \num{1.7} \\
        \bottomrule
    \end{tabular}
    \small
    \label{tab:earth_mars_L-SODA_transcription}
\end{table}
Both approaches maintain failure rates $\beta_{\rm R}$ below the target $\beta$. However, the spectral radius method proves overly conservative: its nominal error is significantly negative, indicating excessive safety margins in the imposed constraints. This leads to an additional fuel consumption of \SI{0.7}{\kilo\gram} compared to the first-order method. This corresponds to a non-negligible penalty given the small magnitude of the uncertainties.
The reported failure rate of the spectral radius variant is $<$\num{3e-3}, since no failed trajectory was observed among the \num{e5} Monte Carlo samples. As discussed by \citet{HanleyLippmanHand_1983_INGWIEAR}, this does not imply a zero failure probability, but only enables the computation of an upper bound on the actual risk.
Furthermore, due to repeated eigenvalue computations on ill-conditioned matrices, the spectral radius variant requires \SI{51}{\%} more runtime than the first-order method for a similar number of iterations. In contrast, first-order variant matches the deterministic solution presented in \citet{CalebEtAl_2025_TPBCSfFOLTTO} in terms of cost, while achieving a low conservatism (\num{1.7}) corresponding to a failure rate of \SI{2.93}{\%}. 

\subsubsection{Failure risk estimations for high-dimensional Gaussian distribution}
\begin{figure}[h]
    \centering 
    \includegraphics[width=0.5\linewidth]{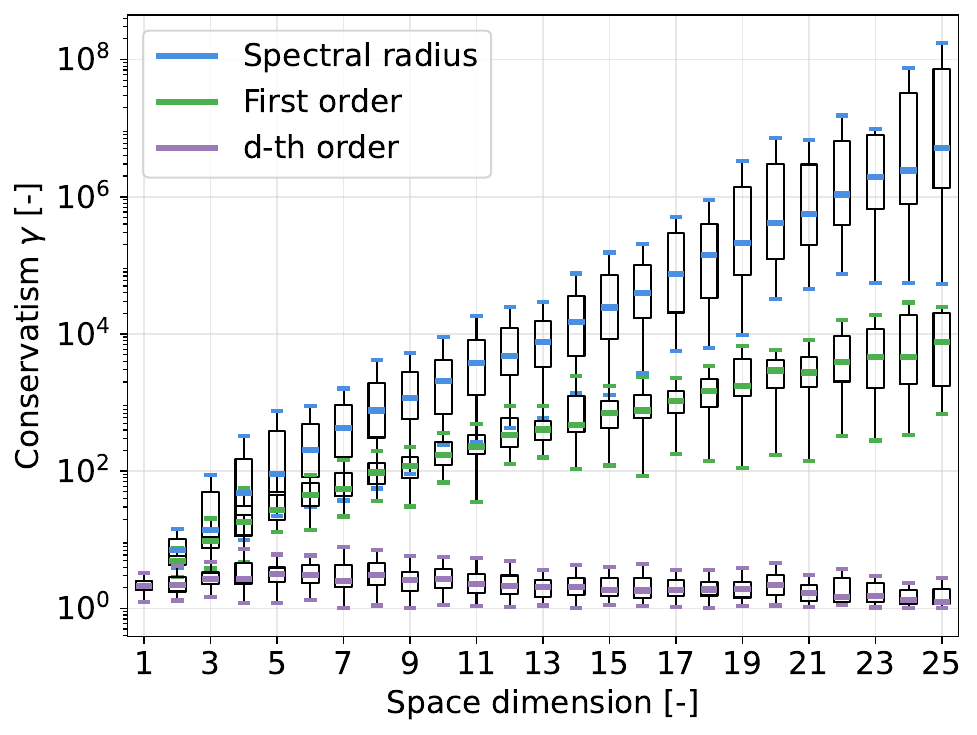}
    \caption{Evolution of the conservatism with the dimension.}
    \label{fig:transcription_high_dimension}
\end{figure}

The three multi-dimensional risk estimation methods developed in this work are now compared using Gaussian distributions of dimension $d$ ranging from \num{1} to \num{25}. No comparison with state-of-the-art methods is possible, as no existing method, to the best of the authors' knowledge, can tackle multi-dimensional Gaussian variables of such size.
For each dimension $d$, $n = \num{1000}$ distributions $\bm{y}\sim\mathcal{N}\left(\Bar{\bm{y}}, \bm{M}\bm{M}^{\textrm{T}}\right)$ are generated randomly. Each component of $\Bar{\bm{y}}$ is drawn from $\mathcal{N}(-1, 0.1)$, and each component of the lower-triangular matrix $\bm{M}$ follows $\mathcal{N}(0, \sigma_m^2)$, with $\sigma_m= \dfrac{\|\Bar{\bm{y}}\|_1}{d^{3/2}\Psi_d^{-1}\left(\beta\right)}$, and $\beta=10^{-3}$.
This scaling ensures that the reference failure probability $\beta_{\textrm{R}}$ remains approximately of the same order of magnitude ($\approx \beta$) regardless of $d$.
Computing $\beta_{\textrm{T},\rho}$, $\beta_{\textrm{T},1}$, $\beta_{\textrm{T},d}$, and $\beta_{\textrm{R}}$ for each distribution yields their associated conservatism values.
\Fig{fig:transcription_high_dimension} presents these results in a boxplot. The central line represents the median, the top and bottom edges the third and first quartiles respectively, and the whiskers extend to \num{1.5} times the interquartile range.
The results show that the spectral radius method is the most conservative, followed by the first-order and then the $d$\textsuperscript{th}-order method, as established in \Thm{thm:conservatism_hierarchy}. The $d$\textendash{th}-order method is the least conservative, especially at high dimension, where it remains stable with conservatism below \num{10}. The other two methods tend to overestimate failure probability near \num{100}\%, causing their conservatism to diverge in high dimensions.

%% file: sections/5_conclusion.tex
\section{Conclusion}
\label{sec:conclusion}

This work addresses the challenge of transcribing multi-dimensional Gaussian chance constraints into tractable deterministic formulations with low conservatism.
We present two transcription methods of increasing efficiency: the spectral radius method, which generalises existing norm-constraint theory to arbitrary multi-dimensional constraints, and the first-order method, which achieves superior tightness with linear computational complexity $\mathcal{O}(d)$.
A $d$-th order risk estimation methodology provides conservative but reliable bounds on the failure probability, applicable to both transcription approaches. This enables \textit{a posteriori} assessment of solution quality without resorting to sampling.

Applied to a standard Earth--Mars transfer optimisation problem with uncertainties, the first-order transcription achieves fuel consumption comparable to deterministic solutions while maintaining the required failure risk $\beta_{\mathrm{R}} = 2.93\%$, below the target $\beta = 5\%$, thereby demonstrating its practical effectiveness.
Conversely, spectral radius methods incur approximately \SI{700}{\gram} of additional fuel consumption due to excessive conservatism, alongside increased computational time arising from the cubic complexity of the transcription.

The proposed $d$-th order risk estimation method provides conservative risk estimates while maintaining bounded conservatism (below $\num{10}$), even in high dimensions. In contrast, previously developed methods exhibit exponential growth in conservatism with respect to constraint dimension. Finally, this risk estimation framework enables a posteriori verification of chance constraint satisfaction without Monte-Carlo simulation, thanks to its quadratic computational complexity.

%% file: sections/7_funding_acknoledgments.tex
\section*{Funding sources}
This work was funded by SaCLaB (grant number 2022-CIF-R-1), a research group of ISAE-SUPAERO. 
\section*{Acknowledgments}
The authors would like to thank Jean-Louis Pac for his suggestions and the careful review of the proofs. 
The authors acknowledge the use of an AI tool (Claude Haiku 4.5) to assist in the reformulation of certain sentences.
\section*{Declaration of competing interest}
The authors have no competing interests to declare that are relevant to the content of this article.

%% file: sections/6_appendix.tex
\section{Appendix}
\label{sec:appendix}

\textbf{Proof of \Thm{thm:spectral_radius_transcription}.} Let $\mathcal{E}\left(1-\beta\right)$ be the iso-\gls*{PDF} ellipsoid of $\bm{y}$ centered in $\Bar{\bm{y}}$ such that $\mathbb{P}\left(\bm{y}\in \mathcal{E}\left(1-\beta\right)\right)=1-\beta$. 
Then, $\bm{v}\in\mathcal{E}\left(1-\beta\right)$ is equivalent to satisfying $\left(\bm{v}-\Bar{\bm{y}}\right)^{\textrm{T}}\bm{\Sigma}_{y}^{-1}\left(\bm{v}-\Bar{\bm{y}}\right)\leq\Phi_d^{-1}\left(1-\beta\right)$.

Consequently, the maximum stretch (semi-major axis) of $\mathcal{E}\left(1-\beta\right)$ is $\sqrt{\Phi_d^{-1}\left(1-\beta\right)}\rho\left(\bm{\Sigma}_{y}\right)$ \citep{Tong_1990_TMND}.
Therefore, we have: $\mathcal{E}\left(1-\beta\right)\subset \mathcal{B}\left(\Bar{\bm{y}}, \Psi_d^{-1}\left(\beta\right)\rho\left(\bm{\Sigma}_{y}\right)\right)$,
where $\mathcal{B}\left(\Bar{\bm{y}}, \Psi_d^{-1}\left(\beta\right)\rho\left(\bm{\Sigma}_{y}\right)\right)$ is the ball centered in $\Bar{\bm{y}}$ of radius $\Psi_d^{-1}\left(\beta\right)\rho\left(\bm{\Sigma}_{y}\right)$. 
Thus: 
\begin{equation*}
    \mathcal{B}\left(\Bar{\bm{y}}, \Psi_d^{-1}\left(\beta\right)\rho\left(\bm{\Sigma}_{y}\right)\right)\subset \left(\mathbb{R}^{-}\right)^d \implies \mathbb{P}(\bm{y}\preceq \bm{0}) \geq  \mathbb{P}(\bm{y}\in \mathcal{E}\left(1-\beta\right)) = 1-\beta.
\end{equation*}
Notice that: $\forall i \in\mathbb{N}_d^*, \ \Bar{y}_i + \Psi_d^{-1}\left(\beta\right)\rho\left(\bm{\Sigma}_{y}\right)  \leq 0 \Leftrightarrow \mathcal{B}\left(\Bar{\bm{y}}, \Psi_d^{-1}\left(\beta\right)\rho\left(\bm{\Sigma}_{y}\right)\right)\subset \left(\mathbb{R}^{-}\right)^d$.

\textbf{Proof of \Rem{rmk:oguri}.}
Consider $y = \|\bm{u}\| - u_{\max}$. Under the weakly non-linear assumption, the norm can be linearised around the nominal control. Applying \Thm{thm:spectral_radius_transcription} yields: 
\begin{equation*}
    \|\Bar{\bm{u}}\| + \Psi_1^{-1}(\beta) \dfrac{\sqrt{\Bar{\bm{u}}^{\textrm{T}}\bm{\Sigma}_{\bm{u}}\Bar{\bm{u}}}}{\|\Bar{\bm{u}}\|} \leq u_{\max} \implies \mathbb{P}\left(\|\Bar{\bm{u}}\|\leq u_{\max}\right) \geq 1 - \beta.
\end{equation*}
We aim to show that the safety bound of \Thm{thm:spectral_radius_transcription} is smaller than that shown in \citet{OguriLantoine_2022_SSCPfRLTTDuU}, in other words: $\Psi_1^{-1}(\beta) \frac{\sqrt{\Bar{\bm{u}}^{\textrm{T}}\bm{\Sigma}_{\bm{u}}\Bar{\bm{u}}}}{\|\Bar{\bm{u}}\|} \leq \Psi_{N_u}^{-1}(\beta) \rho\left(\bm{\Sigma}_{\bm{u}}\right)$.
Indeed, we first note that $\dfrac{\sqrt{\Bar{\bm{u}}^{\textrm{T}}\bm{\Sigma}_{\bm{u}}\Bar{\bm{u}}}}{\|\Bar{\bm{u}}\|}\leq \rho\left(\bm{\Sigma}_{\bm{u}}\right)$, and that the function $\Psi_d^{-1}(\beta)$ is increasing with $d$. 

\textbf{Proof of \Thm{thm:first_order_transcription}.}
Let $\bm{z}\sim\mathcal{N}\left(\bm{0}, \mathbb{I}_d\right)$, since $\bm{\Sigma}_{y}$ is positive definite there is $\bm{M}$ such that $\bm{\Sigma}_{y} = \bm{M}\bm{M}^{\textrm{T}}$ \citep{Cholesky_1910_SLRNDSDL}. Then, $\bm{y} = \Bar{\bm{y}} + \bm{M} \bm{z}$ \citep{Tong_1990_TMND}. Thus:
\begin{equation*}
    \mathbb{P}(\bm{y} \preceq \bm{0}) = \mathbb{P}\left(\bm{M} \bm{z} + \Bar{\bm{y}} \preceq \bm{0}\right) = \mathbb{P}\left(\bigcap_{i=1}^{d}{\left\{\bm{M}_i \bm{z} + \Bar{y}_i\leq 0 \right\}}\right),
\end{equation*}
where $\bm{M}_i$ is the $i$\textendash{th} row of $\bm{M}$.
Note that the constraints in the Cartesian space: $\bm{y} \preceq\bm{0}$ can be expressed using the following transformation $\forall \bm{v}\in\mathbb{R}^d$: $\bm{v}' = \bm{M}^{-1} \left(\bm{v}-\Bar{\bm{y}}\right)$, inspired by the Mahalanobis distance \citep{Mahalanobis_1936_OtGDiS}.
Indeed, the Euclidean norm of a transformed point $\bm{v}'$ is the Mahalanobis distance between $\bm{v}$ and the distribution $\bm{y}$: $\left\|\bm{v}'\right\|^2 = \bm{v}'^{\textrm{T}} \bm{v}'=\left(\bm{v}-\Bar{\bm{y}}\right)^{\textrm{T}} \bm{\Sigma}_{y}^{-1} \left(\bm{v}-\Bar{\bm{y}}\right)$. In addition it transforms $\Bar{\bm{y}}$ into $\bm{0}$, and $\bm{y}$ into $\bm{y}'$, which follows the same distribution as $\bm{z}$.

Let $i\in\mathbb{N}_d^* $, each equation $\bm{M}_i\bm{v}' + \Bar{y}_i = 0$ of unknown $\bm{v}'$ is the equation of an hyper-plane $\mathcal{H}_i'$ of $\mathbb{R}^d$, and the vector $\bm{M}_{i}$ is normal to $\mathcal{H}_i'$.
Therefore the Euclidean distance between any point $\bm{v}'$ and $\mathcal{H}_i'$ is: $ d\left(\bm{v}', \mathcal{H}_i'\right) = \dfrac{\left\vert \bm{M}_i\bm{v}' + \Bar{y}_i\right\vert }{\|\bm{M}_{i}\|}$.
Consequently: $ r_i = d\left(\Bar{\bm{y}}', \mathcal{H}_i'\right) = \dfrac{\vert \Bar{y}_i\vert }{\|\bm{M}_i\|}$.
The value $r_i$ is the distance between $\Bar{\bm{y}}'$ and $\mathcal{H}'_i$.
Note that $\mathcal{H}_i$, the antecedent of the hyper-plane $\mathcal{H}_i'$, corresponds to the hyper-plane of equation $v_i=0$, thus, to the $i$\textendash{th} edge of the admissible domain.
In other words, the Mahalanobis distance between $\Bar{\bm{y}}$ and the $i$\textendash{th} constraint $\mathcal{H}_i$ in the Cartesian space is $r_i$.

If $r_i \geq \Psi_d^{-1}\left(\beta\right)$, then $\mathcal{B}\left(\Bar{\bm{y}}', \Psi_d^{-1}\left(\beta\right)\right)$, the ball centered in $\Bar{\bm{y}}'$ of radius $\Psi_d^{-1}\left(\beta\right)$, never intersects with $\mathcal{H}'_i$ or is tangent. 
Thus, if $\Bar{\bm{y}} \preceq \bm{0}$ and $\forall i\in\mathbb{N}_d^*, \ r_i \geq\Psi_d^{-1}\left(\beta\right)$:
\begin{equation*}
    \mathcal{B}\left(\Bar{\bm{y}}', \Psi_d^{-1}\left(\beta\right)\right) \subset \bigcap_{i=1}^{d}{\left\{\bm{v}'\in\mathbb{R}^d \ | \ \bm{M}_i\bm{v}' + \Bar{y}_i \leq 0 \right\}}.
\end{equation*}
Since $\|\bm{z}\|^2$ follows a $\chi^2$ distribution with $d$ degrees of freedom \citep{Tong_1990_TMND}, if $\Bar{\bm{y}} \preceq \bm{0}$ and $\forall i\in\mathbb{N}_d^*, \ r_i \geq\Psi_d^{-1}\left(\beta\right)$, it leads to:
\begin{equation*}   
         \mathbb{P}(\bm{y} \preceq \bm{0})=\mathbb{P}\left(\bigcap_{i=1}^{d}{\left\{\bm{M}_i\bm{y}' + \Bar{y}_i \leq 0 \right\}}\right) \geq  \mathbb{P}\left(\bm{y}'\in\mathcal{B}\left(\Bar{\bm{y}}', \Psi_d^{-1}\left(\beta\right)\right)\right) = \mathbb{P}\left(\|\bm{z}\|^2 \leq \Phi_d^{-1}\left(1-\beta\right)\right) =  1-\beta.
\end{equation*}
Therefore: $\Bar{\bm{y}} \preceq \bm{0} \land \left(\forall i\in\mathbb{N}_d^*, \ \frac{\vert \Bar{y}_i\vert }{\|\bm{M}_i\|} \geq\Psi_d^{-1}\left(\beta\right)\right) \implies \mathbb{P}(\bm{y} \preceq \bm{0}) \geq 1 - \beta$.
To obtain the result, note that $\sigma_{y,i}^2=\Sigma_{\bm{y},i,i}=\bm{M}_i\bm{M}_i^{\textrm{T}}=\|\bm{M}_i\|^2$, and: 
\begin{equation*}
    \forall i\in\mathbb{N}_d^*, \ \Bar{\bm{y}}_i  + \Psi_d^{-1}\left(\beta\right)\bm{\sigma}_{y}  \leq \bm{0} \Leftrightarrow \Bar{\bm{y}} \preceq \bm{0} \land \left(\forall i\in\mathbb{N}_d^*, \ \frac{\vert \Bar{y}_i\vert }{\|\bm{M}_i\|} \geq\Psi_d^{-1}\left(\beta\right)\right).
\end{equation*}

\textbf{Proof of \Rem{rmk:spectral_vs_first_order}.} 
If $d=1$, then $\sigma_{y} = \sqrt{\Sigma_y}=\rho(\Sigma_y)$.
If $d > 1$, it suffices to prove that the margin introduced by the spectral radius transcription is greater than that of the first-order transcription: for all $i$, $\sigma_{\bm{y},i} \leq \rho(\bm{\Sigma}_{\bm{y}})$.  
Indeed, if $\bm{v}\in\mathbb{R}^d$ such that $\|\bm{v}\|=1$, we have: $\sqrt{\bm{v}^{\textrm{T}}\bm{\Sigma}_{\bm{y}}\bm{v}}\leq\rho\left(\bm{\Sigma}_{\bm{y}}\right)$.
Taking $\bm{v}$ to be the $i$-th canonical basis vector yields the result.

\textbf{Proof of \Prop{prop:risk_estimates_all}.}
Each estimate results from a deterministic sufficient condition of the form $\Bar{\bm{y}} + \bm{\Delta} \bm{y}(\Bar{\bm{y}},\bm{\Sigma}_{\bm{y}},\beta_{\textrm{T}}) \preceq \bm{0}$, for some deterministic safety bound $\bm{\Delta}\bm{y}(\Bar{\bm{y}},\bm{\Sigma}_{\bm{y}},\beta_{\textrm{T}})$ increasing in $\beta_{\textrm{T}}$. In each case, if the condition holds for a given $\beta_{\textrm{T}}$, then by construction, we have: $\mathbb{P}(\bm{y} \preceq \bm{0}) \geq 1 - \beta_{\textrm{T}} \implies \beta_{\textrm{T}} \geq \beta_{\textrm{R}}$.
For instance, for the linear univariate Gaussian risk estimation from \citet{NakkaChung_2023_TOoCCNSSfMPuU}, $\beta_{\textrm{T},\textrm{NC}}$ is obtained by solving the limit case of the sufficient condition $\bm{h}^{\textrm{T}}\Bar{\bm{z}} + \sqrt{\frac{1-\beta_{\textrm{T}}}{\beta_{\textrm{T}}}}\sqrt{\bm{h}^{\textrm{T}}\bm{\Sigma}_{\bm{z}}\bm{h}} = a$ for $\beta_{\textrm{T}}$. All other risk estimations were derived similarly. Note that since \citet{BlackmoreEtAl_2011_CCOPPwO} provide a necessary and sufficient condition, $\beta_{\textrm{T},\textrm{B}}=\beta_{\textrm{R}}$.

\textbf{Proof of \Thm{thm:dth_order_risk}.}
We first introduce the following definitions.
\begin{enumerate}[label=(\alph*), itemjoin={,\ }]
    \item If $R\geq 0$ and $\bm{z}\in\mathbb{R}^d$, $\mathcal{B}\left(\bm{z}, R\right)$ is the ball centred in $\bm{z}$ of radius $R$.
    \item If $R_1\leq R_2$, $\Delta\mathcal{B}\left(\bm{z}, R_1, R_2\right) = \mathcal{B}\left(\bm{z}, R_2\right) \smallsetminus \mathcal{B}\left(\bm{z}, R_1\right)$.
    \item $S\left(\bm{z}, R, \bm{v}, \theta\right)$ is a spherical sector of $\mathcal{B}\left(\bm{z}, R\right)$ of angle $\theta$ and axis $\bm{v}\in\mathbb{R}^d$.
    \item $\Delta S\left(\bm{z}, R_2, R_1, \bm{v}, \theta\right)=S\left(\bm{z}, R_2, \bm{v}, \theta\right)\smallsetminus S\left(\bm{z}, R_1, \bm{v}, \theta\right)$.
    \item $\Delta\Psi_d\left(R_{1},R_{2}\right) = \Psi_d\left(R_{1}\right) - \Psi_d\left(R_{2}\right)$. If $0\leq j<i\leq d$, the spherical sector $S_{i,j}=S\left(\Bar{\bm{y}}', \Tilde{r}_i, \Tilde{\bm{M}}_j, \theta_{i,j}\right)$ is such that $\mathcal{B}\left(\Bar{\bm{y}}', \Tilde{r}_i\right)\smallsetminus S_{i,j}$ does not cross the constraint $\Tilde{\mathcal{H}}_j'$.
    \item For short: $\mathcal{I}\left(x\right)=I_{\sqrt{1-x}}\left(\frac{d-1}{2}, \frac{1}{2}\right)$.
\end{enumerate}
We first prove that if $\bm{z}\sim\mathcal{N}\left(\bm{0}, \mathbb{I}_d\right)$, $R_2\geq R_1\geq0$, $\bm{v}\in\mathbb{R}^d$, and $\theta\in[0, \frac{\pi}{2}]$, then: 
\begin{equation*}
    \mathbb{P}\left(\bm{z}\in \Delta S\left(\bm{0}, R_2, R_1, \bm{v}, \theta\right)\right) = \frac{\mathcal{I}\left( \cos\theta\right) }{2}\Delta\Psi_d\left(R_{1},R_{2}\right).
\end{equation*}
Indeed, the volume of a spherical sector $S(\bm{0}, R, \bm{v}, \theta)$ is $\mathcal{S}\left(R, \theta\right) = \frac{1}{2} \mathcal{V}\left(R\right)\mathcal{I}\left( \cos\theta\right)$ \citep{Li_2011_CFftAaVoaHC}, where $\mathcal{V}\left( R\right)$ is the volume of $\mathcal{B}\left(\bm{x}, R\right)$ \citep{LejeuneDirichlet_1839_SUNMPLDDIM}.
The probability for $\bm{z}\sim\mathcal{N}\left(\bm{0},\mathbb{I}_d\right)$ to belong to a sector $S\left(\bm{0}, R, \bm{v}, \theta\right)$ knowing that $\bm{z}\in\mathcal{B}\left(\bm{0}, R\right)$ follows a uniform distribution of size $\mathcal{V}\left(R\right)$,
Therefore: $ \mathbb{P}\left(\bm{z}\in S(\bm{0}, R, \bm{v}, \theta) | \ \bm{z}\in\mathcal{B}\left(\bm{0}, R\right)\right) = \frac{\mathcal{S}\left(R, \theta\right)}{\mathcal{V}\left(R\right)}=\frac{1}{2} \mathcal{I}\left( \cos\theta\right)$.
Since $S(\bm{0}, R, \bm{v}, \theta)\subset \mathcal{B}(\bm{0}, R)$:
\begin{equation*}
    \mathbb{P}\left(\bm{z}\in S(\bm{0}, R, \bm{v}, \theta)\right) = \mathbb{P}\left(\bm{z}\in \mathcal{B}\left(\bm{0}, R\right)\right)  \mathbb{P}\left(\bm{z}\in S(\bm{0}, R, \bm{v}, \theta) | \ \bm{z}\in\mathcal{B}\left(\bm{0}, R\right)\right) = \frac{1}{2}\Phi_{d}\left(R^2\right)\mathcal{I}\left( \cos\theta\right).
\end{equation*}
Let $R_2\geq R_1\geq0$, and two sectors $S(\bm{0}, R_1, \bm{v}, \theta)\subset S(\bm{0}, R_2, \bm{v}, \theta)$, since $\Phi_{d}\left(R^2\right)=1-\Psi_d\left(R\right)$, then: 
\begin{equation*}
    \begin{aligned}
         \mathbb{P}\left(\bm{z}\in \Delta S\left(\bm{0}, R_2, R_1, \bm{v}, \theta\right)\right) = & \mathbb{P}\left(\bm{z}\in S(\bm{0}, R_2, \bm{v}, \theta)\right) -\mathbb{P}\left(\bm{z}\in S(\bm{0}, R_1, \bm{v}, \theta)\right) \\ = & \frac{\mathcal{I}\left( \cos\theta\right) }{2}\Delta\Psi_d\left(R_{2}, R_{1}\right).
    \end{aligned}
\end{equation*}
This concludes our intermediate proof. 

The rest of this proof uses the objects $\bm{y}'$, $\bm{M}$, and $\mathcal{H}'_{i}$ from the proof of \Thm{thm:first_order_transcription} where more details can be found. We define the function $\alpha$ as the permutation function such that $\forall i\in\mathbb{N}_d, \ r_{\alpha(i)} = \Tilde{r}_i$, by convention $\Tilde{r}_0 = 0$, we also denote $\Tilde{\mathcal{H}}'_i=\mathcal{H}'_{\alpha(i)}$ and $\Tilde{\bm{M}}_i=\bm{M}_{\alpha(i)}$ respectively the edge of the $i$\textendash{th} closest constraint to $\Bar{\bm{y}}'$ in the transformed space and one of its normal vectors.
For $d=1$: 
\begin{equation*}
    \beta_{\textrm{T},d} =   1 - \left[\Psi_{d}\left(\Tilde{r}_{0}\right)-\Psi_{d}\left(\Tilde{r}_{1}\right)\right] \max\left[0, 1 - 0 \right]  = \Psi_{d}\left(\frac{-\Bar{y}}{\sqrt{\Sigma}_y}\right) = \beta_{\textrm{T},1} \geq \beta_{\textrm{R}}.
\end{equation*}
Let us now focus on $d>1$.       
Let $i\in\mathbb{N}_d^*$, the proof of \Thm{thm:first_order_transcription} shows that $\mathcal{B}\left(\bm{0}, \Tilde{r}_i\right)$ intersects with $i-1$ hyper-planes: $\Tilde{\mathcal{H}}'_{1}, \cdots, \Tilde{\mathcal{H}}'_{i-1}$. Each of these hyper-planes is at a distance $\Tilde{r}_1, ..., \Tilde{r}_{i-1}\leq\Tilde{r}_i$ from $\bm{0}$. Therefore, for each intersecting hyper-plane $\Tilde{\mathcal{H}}'_{j}$, with $1\leq j<i$, the sector $S_{i,j}$, with $\theta_{i,j}=\arccos \frac{\Tilde{r}_j}{\Tilde{r}_i}$, does not intersect with $\Tilde{\mathcal{H}}'_{j}$.    
In other words, after slicing the sector of $\mathcal{B}\left(\bm{0}, \Tilde{r}_i\right)$ that has its spherical cap on the wrong side of $\Tilde{\mathcal{H}}'_{j}$, the remaining set satisfies the constraints. Thus, if $\Omega'$ is the transformed admissible set and $\Bar{\bm{y}}\preceq \bm{0}$: 
\begin{equation*}
    \forall i\in\mathbb{N}_d^*, \  \mathcal{B}\left(\bm{0}, \Tilde{r}_i\right)  \smallsetminus \bigcup_{j=1}^{i-1} {S_{i,j}}\subset \Omega' \implies \bigcup_{i=1}^{d}{\left(\mathcal{B}\left(\bm{0}, \Tilde{r}_i\right)\smallsetminus \bigcup_{j=1}^{i-1} {S_{i,j}}\right)} \subset \Omega'.
\end{equation*}
Note that if $\Delta S_{i,j} = S_{i,j}\cap\Delta\mathcal{B}\left(\bm{0}, \Tilde{r}_i, \Tilde{r}_{i-1}\right)$:
\begin{equation*}
    \bigcup_{i=1}^{d}{\left(\mathcal{B}\left(\bm{0}, \Tilde{r}_i\right)\smallsetminus \bigcup_{j=1}^{i-1} {S_{i,j}}\right)} = \bigsqcup_{i=1}^{d}{\left(\Delta\mathcal{B}\left(\bm{0}, \Tilde{r}_i, \Tilde{r}_{i-1}\right)\smallsetminus \bigcup_{j=1}^{i-1} {\Delta S_{i,j}}\right)}.
\end{equation*}
However, in general $\bigcap_{j=1}^{i-1} {\Delta S_{i,j}} \neq \emptyset$, therefore: $\sum_{j=1}^{i-1}{\mathbb{P}\left(\bm{z}\in\Delta S_{i,j}\right)\geq\mathbb{P}\left(\bm{z}\in\bigcup_{j=1}^{i-1} {\Delta S_{i,j}}\right)}$. Consequently, since $\Bar{\bm{y}}\preceq \bm{0}$: 
\begin{equation*}
    \begin{aligned}
        \mathbb{P}\left(\bm{z}\in\Omega'\right) \geq &  
        \mathbb{P}\left(\bm{z}\in\bigsqcup_{i=1}^{d}{\left(\Delta\mathcal{B}\left(\bm{0}, \Tilde{r}_i, \Tilde{r}_{i-1}\right)\smallsetminus \bigcup_{j=1}^{i-1} {\Delta S_{i,j}}\right)}\right)
         \geq \sum_{i=1}^{d}{\mathbb{P}\left(\bm{z}\in{\left(\Delta\mathcal{B}\left(\bm{0}, \Tilde{r}_i, \Tilde{r}_{i-1}\right)\smallsetminus \bigcup_{j=1}^{i-1} {\Delta S_{i,j}}\right)}\right)} \\
        \geq & \sum_{i=1}^{d}\left[\mathbb{P}\left(\bm{z}\in{\Delta\mathcal{B}\left(\bm{0}, \Tilde{r}_i, \Tilde{r}_{i-1}\right)}\right) \mathbb{P}\left(\bm{z}\in\bigcup_{j=1}^{i-1} {\Delta S_{i,j}}\right)\right] \\ \geq & \sum_{i=1}^{d}\max\Bigg[0, \mathbb{P}\left(\bm{z}\in{\Delta\mathcal{B}\left(\bm{0}, \Tilde{r}_i, \Tilde{r}_{i-1}\right)}\right) -  \sum_{j=1}^{i-1}{\mathbb{P}\left(\bm{z}\in\Delta S_{i,j}\right)} \Bigg] \\
         \geq & \sum_{i=1}^{d}\max\Bigg[0, \Delta\Psi_{d}\left(\Tilde{r}_{i-1},\Tilde{r}_{i}\right)  - \sum_{j=1}^{i-1}{\frac{\mathcal{I}\left( \cos\theta_{i,j}\right) }{2}\Delta\Psi_{d}\left(\Tilde{r}_{i-1},\Tilde{r}_{i}\right)} \Bigg] =  1-\beta_{\textrm{T},d} \\
    \end{aligned}
\end{equation*}
Since $\Delta\Psi_{d}\left(\Tilde{r}_{i-1},\Tilde{r}_{i}\right) \geq 0$, and $\bm{y}'$ and $\bm{z}$ follow the same distribution: $1-\beta_{\textrm{R}} =\mathbb{P}\left(\bm{y}'\in\Omega'\right)
        \geq  1-\beta_{\textrm{T},d}$

\textbf{Proof of \Thm{thm:conservatism_hierarchy}.}
If $d=1$, the first-order and spectral radius transcriptions coincide, hence $\gamma_1 = \gamma_\rho$. In addition, the proof of \Thm{thm:dth_order_risk} in \Thm{thm:dth_order_risk} shows that $\beta_{\mathrm{T},d} = \beta_{\mathrm{T},1}$, so that $\gamma_d = \gamma_1$.
When $d > 1$, we establish the inequality $\beta_{\mathrm{R}} \leq \beta_{\mathrm{T},d} \leq \beta_{\mathrm{T},1} \leq \beta_{\mathrm{T},\rho}$. Since $\gamma$ is monotonically increasing, this directly yields the result.
The inequality $\beta_{\mathrm{T},1} \leq \beta_{\mathrm{T},\rho}$ follows from the fact that $\sigma_{\bm{y},i} \leq \rho(\bm{\Sigma}_{\bm{y}})$ for all $i \in \{1,\dots,d\}$.
The bound $\beta_{\mathrm{R}} \leq \beta_{\mathrm{T},d}$ is guaranteed by \Thm{thm:dth_order_risk}.
Finally, the difference between the first-order and $d$th-order estimates is given by: 
\begin{equation*}
    \beta_{\textrm{T},1} - \beta_{\textrm{T},d} = \sum_{i=2}^d \left[\Psi_d(\Tilde{r}_i) - \Psi_d(\Tilde{r}_{i-1})\right] \max\left[0, 1 - \frac{1}{2} \sum_{j=1}^{i-1} I\left(\frac{d-1}{2}, \frac{1}{2}; \frac{\Tilde{r}_j}{\Tilde{r}_i}\right)\right] \geq 0.
\end{equation*}